\documentclass{amsart}
\usepackage{palatino}
\usepackage{amsthm, amsmath}
\usepackage{amssymb} 
\usepackage{amscd}
\usepackage[margin=3.5cm]{geometry}
\usepackage[breaklinks,bookmarksopen,bookmarksnumbered]{hyperref}
\usepackage[all]{xy}

\theoremstyle{plain}

\newtheorem{lem}{Lemma}
\newtheorem{prop}{Proposition}
\newtheorem*{prop*}{Proposition}
\newtheorem*{cor*}{Corollary}
\theoremstyle{remark}
\newtheorem{rem}{Remark}
\theoremstyle{remark}

\newcommand\pr{\noindent\textit{Proof} : }

\newcommand\rond{\kern 1pt{\scriptstyle\circ}\kern 1pt}

\newcommand\Aut{\operatorname{Aut}}

\newcommand\Ker{\operatorname{Ker}}

\newcommand\rk{\operatorname{rk}}
\newcommand\Tr{\operatorname{Tr}}

\newcommand\Z{\mathbb{Z}}
\newcommand\Q{\mathbb{Q}}
\newcommand\R{\mathbb{R}}
\newcommand\C{\mathbb{C}}
\renewcommand\P{\mathbb{P}}

\newcommand\F{\mathbb{F}}

\renewcommand\O{\mathcal{O}}
\renewcommand\Q{\mathcal{Q}}

\newcommand\iso{\vbox{\hbox to .8cm{\hfill{$\scriptstyle\sim$}\hfill}
\nointerlineskip\hbox to .8cm{{\hfill$\longrightarrow $\hfill}} }}

\begin{document}
\baselineskip15pt
\title[Abelian varieties associated to  Gaussian lattices]{Abelian varieties associated to  Gaussian lattices}
\author[Arnaud Beauville]{Arnaud Beauville}
\address{Laboratoire J.-A. Dieudonn\'e\\
UMR 7351 du CNRS\\
Universit\'e de Nice\\
Parc Valrose\\
F-06108 Nice cedex 2, France}
\email{arnaud.beauville@unice.fr}
 
\date{\today}

\begin{abstract}
We associate to a unimodular lattice $\Gamma $, endowed with an automorphism of square $-1$, a principally polarized abelian variety $A_\Gamma =\Gamma _\R/\Gamma $. We show that the configuration of $i$-invariant theta divisors of $A_\Gamma $ follows a pattern very similar to the classical theory of theta characteristics; as a consequence we find that $A_\Gamma $ has a high number of vanishing theta\-nulls. When $\Gamma =E_8$ we recover the 10 vanishing thetanulls of the abelian fourfold discovered by R.~Varley.
\end{abstract}
\maketitle 
\section*{Introduction}
A \emph{Gaussian lattice} is a free, finitely generated  $\Z[i]$-module $\Gamma $ with a positive hermitian form $\Gamma \times \Gamma \rightarrow \Z[i]$. Equivalently, we can view $\Gamma $ as a lattice over $\Z$ endowed with an automorphism $i$ of square $-1_{\Gamma }$. This gives a complex structure on the vector space $\Gamma _\R:=\Gamma \otimes_{\Z} \R$; we associate to   $\Gamma $ the complex torus $A_\Gamma :=\Gamma_{\R}/\Gamma $. 
\par As a complex torus $A_\Gamma $ is isomorphic to $E^g$, where $E$ is the complex elliptic curve $\C/\Z[i]$ and $g=\frac{1}{2} \rk_{\Z}\Gamma $. More interestingly,
the hermitian form provides a \emph{polarization} on $A_\Gamma $ (see (\ref{abvar}) below); in particular, if $\Gamma $ is unimodular, $A$ is a principally polarized abelian variety (p.p.a.v. for short), which is indecomposable if $\Gamma $ is indecomposable.
\par The first non-trivial case is $g=4$, with $\Gamma $ the root lattice of type $E_8$ (Example \ref{examples}.1).
The resulting p.p.a.v. is the abelian fourfold discovered by Varley \cite{V} with a different (and more geometric) description; it has 10 ``vanishing thetanulls'' (even theta functions vanishing at 0),
the maximum possible for a 4-dimensional indecomposable p.p.a.v. In fact this property characterizes the Varley fourfold outside the hyperelliptic Jacobian locus \cite{D}.
\par Our aim is to explain this property from the lattice point of view, and  to extend it  to all unimodular lattices. It turns out that  we can mimic  the classical theory of theta characteristics, replacing the automorphism $(-1)$ by $i$. We will show:

\par $\bullet$ The group $A_i$ of $i$-invariant  points of $A_\Gamma $ is a $\F_2$-vector space of dimension $g$; it admits a natural non-degenerate bilinear symmetric form $b$.
\par $\bullet$ The set  of  $i$-invariant theta divisors of $A_\Gamma $ is an affine space over $A_i$, isomorphic to the space of quadratic forms on $A_i$ associated to $b$ (see (\ref{f2})).
\par $\bullet$ Let $\Theta$ be an $i$-invariant theta divisor, and $Q$  the corresponding quadratic form. The multiplicity $m_0(\Theta )$ of $\Theta $ at $0$ satisfies
\[2m_0 (\Theta ) \equiv \sigma (Q)+g\quad \mathrm{(mod.\ 8)}\, ,\]
where $\sigma $ is the \emph{Brown invariant} of the form $Q$ (\ref{f2}).

\par As a consequence, we obtain a high number of $i$-invariant divisors  $\Theta $ with $m_0(\Theta )\equiv 2$ $(\mathrm{mod.}\ 4)$; each of them corresponds to a 
vanishing thetanull. When $\Gamma $ is even, this number is $2^{\frac{g}{2} -1}(2^{\frac{g}{2} }-(-1)^{\frac{g}{4} })$; for $g=4$ we recover the 10 vanishing thetanulls of the Varley fourfold.

\medskip	
\section{Gaussian lattices}
\subsection{Lattices}  As recalled in the Introduction, a Gaussian lattice is a free finitely generated  \hbox{$\Z[i]$-module}  $\Gamma $ endowed with a positive hermitian form\footnote{Our convention is that $H(x,y)$ is $\C$-linear in $y$.} $H:\Gamma \times \Gamma \rightarrow \Z[i]$. We write $H(x,y)=\allowbreak S(x,y)+iE(x,y)$; $S$ and $E$ are $\Z$-bilinear forms on $\Gamma $, $S$ is symmetric, $E$ is skew-symmetric, and we have
\[ S(ix,iy)= S(x,y) \quad,\quad  E(ix,iy)= E(x,y) \quad,\quad E(x,y)=S(ix,y)\ .\]

We will rather view a Gaussian lattice as an ordinary lattice (over $\Z$) with an automorphism $i$ such that $i^2=-1_\Gamma $: the last formula above defines $E$, and we have  $H=S+iE$.

\par We have $\det S =\det E =(\det H)^2$; the lattice is \emph{unimodular} when these numbers are equal to 1. It is {\it even} if $S(x,x)$ is even for all $x\in \Gamma $.   We say  that $\Gamma $ is \emph{indecomposable over} $\Z[i]$ if it  cannot be written as the orthogonal sum of two nonzero Gaussian lattices; this is of course the case if $\Gamma $ is indecomposable over $\Z$, but the converse is false (Example 3 below). 

\subsection{Examples}\label{examples}

1) For $ g$ even, the lattice $\Gamma _{2g}$  is  
\[\Gamma _{2g}:= \{(x_j)\in \R^{2g}\ |\ x_j\in \frac{1}{2}\Z\  ,\   x_j - x_k \in \Z\  ,\  \sum x_j\in 2\Z \}\ .\]
The inner product is inherited from the euclidean structure of $\R^{2g}$, and the automorphism $i$ is given in the standard basis $(e_j)$ by
\[ie_{2j-1}=e_{2j}\qquad ie_{2j}=-e_{2j-1}\quad\hbox{for }\ 1\leq j\leq g\ . \]
The lattice $\Gamma _{2g}$ is unimodular, indecomposable when $g>2$, and even if $g$ is divisible by 4. The first case $g=4$ gives the root lattice $E_8$.
\par The automorphism $i$ is \emph{unique} up to conjugacy: for $g=4$ this is classical \cite{C}, and for $g\geq 6$ this follows easily from the fact that $\Aut(\Gamma _{2g})$ is the semi-direct product $(\Z/2)^{2g-1}\rtimes \mathfrak{S}_{2g}$, 
acting by permutation and even changes of sign of the basis vectors $(e_j)$.

\medskip	
\par 2) The Leech lattice $\Lambda_{24} $ admits an automorphism of square $-1$ \cite{C-S}, also unique up to conjugacy.

\medskip	
\par 3) Let $\Gamma _0$ be a   lattice, and $\Gamma := \Gamma _0\otimes_{\Z}\Z[i]$. The inner product of $\Gamma _0$ extends to an hermitian  inner product on $\Gamma $, which is then a gaussian lattice. If $\Gamma _0$ is unimodular, resp.\   even, resp.\  indecomposable, $\Gamma $ is unimodular, resp.\   even, resp.\  indecomposable over $\Z[i]$.

\subsection{The abelian variety \boldmath{$A_\Gamma $}}\label{abvar}
Let $\Gamma $ be a Gaussian lattice, of rank $2g$ over $\Z$. We put $\Gamma _{\R}:=\Gamma \otimes _{\Z}\R$ and $A_\Gamma :=\Gamma_\R/\Gamma $.
The automorphism $i$ defines a complex structure on  $\Gamma _{\R}$, so that $A_\Gamma $ is a complex torus. Since $\Gamma $ is a free $\Z[i]$-module, $A_\Gamma $ is isomorphic to $E^g$, where $E$ is the complex elliptic curve $\C/\Z[i]$.

\par The positive hermitian form $H$ extends to $\Gamma _\R$, and its imaginary part $E$ takes integral values on $\Gamma $: this is by definition a \emph{polarization} on $A_\Gamma $. The polarization is principal if and only if $\Gamma $ is unimodular; the p.p.a.v.
 $A_\Gamma $ is indecomposable (i.e. is not a product of two nontrivial p.p.a.v.) if and only if $\Gamma $ is indecomposable over $\Z[i]$.

\par The multiplication by $i$ on $\Gamma _\R$ induces an automorphism of $A_\Gamma $, that we simply denote $i$.
Conversely, let $A=V/\Gamma $ be a complex torus, of dimension $g$, with an automorphism  inducing on $T_0(A)=V$ the multiplication by $i$. Then $\Gamma $ is a  $\Z[i]$-module, thus isomorphic to $\Z[i]^g$, so that $A$ is isomorphic to $E^g$;  polarizations of $A$ correspond bijectively to positive hermitian forms on $\Gamma $. 

\medskip

\section{Linear algebra over $\F_2[i]$}

\subsection{Linear algebra over $\F_2$}\label{f2}

We consider a vector space $V$ over $\F_2$, of dimension $g$, with
a non-degenerate symmetric bilinear form $b$ on  $V$.  Two different situations may occur:

$\bullet$ $b(x,x)=0$ for all $x\in V$; in that case $b$ is a symplectic form.

$\bullet$ $b(x,x)$ is not identically zero; it is then easy (using induction on $g$) to prove that $V$ admits an orthonormal basis with respect to $b$.

A \emph{quadratic form associated to $b$} is a function $q:V\rightarrow \Z/4$ such that
\[q(x+y)=q(x)+q(y)+2 b(x,y)\qquad \mathrm{for}\ x,y\in V\ ,\]
where multiplication by 2  stands for  the isomorphism $\Z/2\iso 2\Z/4\Z\subset \Z/4\Z$.

Observe that this implies $\ q(0)=0\ $ and $\ q(x)\equiv b(x,x)\ (\mathrm{mod.}\ 2)$. We denote by $\mathcal{Q}_b$  the set  of quadratic forms associated to $b$; it   is an affine space over $V$, the action of $V$ on $\mathcal{Q}_b$ being given by $(\alpha +q)(x)=q(x)+2b(\alpha ,x)$ for $q\in \mathcal{Q}_b$, $\alpha ,x\in V$. 

When $b$ is symplectic, $q$ takes it values in $2\Z/4\Z\cong \Z/2$; the corresponding form $q':V\rightarrow \Z/2$ is a quadratic form associated to $b$ in the usual sense, that is satisfies
$\ q'(x+y)=q'(x)+q'(y)+\allowbreak b(x,y)\ $ for  $x,y\in V$.

The \emph{Brown invariant} $\sigma(q)\in \Z/8$ of a form $q\in\mathcal{Q}_b$ has been introduced in \cite{B} as a generalization of the Arf invariant; it can be defined as follows. If $b$ is symplectic,  we put $\sigma (q):=4\,\mathrm{Arf}(q')$, where $q':V\rightarrow \Z/2$ is the form defined above. Otherwise $b$ admits an orthonormal basis $(e_1,\ldots ,e_g)$; we have $q(e_i)=\pm 1$, and we let $g^+$ (resp. $g^-$)  be the number of basis vectors $e_i$ such that $q(e_i)= 1$ (resp. $-1$). Then $\sigma (q)=g^+-g^-\ (\mathrm{mod.}\ 8)$.

\subsection{Linear algebra over \boldmath{$\mathrm{F}_2[i]$}}\label{f2i}

Let $\Gamma $ be a  unimodular Gaussian lattice of rank $2g$ over $\Z$. We put $A_2:=\Gamma /2\Gamma $; this is naturally identified with 
 the 2-torsion subgroup  of $A_\Gamma $. We have the following structures on $A_2$:

\par $a)$ $A_2$ is a free $\F_2[i]$-module of rank $g$. We put $\varepsilon :=1+i$ in  $\F_2[i]$; then $\F_2[i]=\F_2[\varepsilon ]$, with $\varepsilon ^2=0$. The subgroup $A_i$ of $i$-invariant elements is $\Ker \varepsilon =\varepsilon A_2 $; it is a $\F_2$-vector space of dimension $g$.

\par $b)$ The form $E$ induces  on $A_2$ a symplectic form $e$ (the Weil pairing for $A_\Gamma $). Since $E(x,iy)=\allowbreak -E(ix,y)$,  we have, for $\alpha ,\beta \in A_2$,
\[e(\alpha ,\varepsilon \beta )=e(\varepsilon \alpha ,\beta )\quad\hbox{hence}\quad e(\varepsilon \alpha ,\varepsilon \beta )=0\ ; \]
thus $A_i$ is a Lagrangian subspace of $A_2$.

\par $c)$ The form $x\mapsto S(x,x) $ induces a quadratic form $Q:A_2\rightarrow \Z/4$ associated with the bilinear symmetric form $(\alpha ,\beta )\mapsto e(\alpha ,i\beta )$ (\ref{f2}). In particular we have $Q(\alpha )\equiv e(\alpha ,i\alpha )\ (\mathrm{mod.}\ 2)$.

 Since $S((1+i)x,(1+i)x)=2S(x,x)$, we have $Q(\varepsilon \alpha )=2Q(\alpha )=2e(\alpha ,i\alpha )$.

\begin{lem}
Let $q:A_2\rightarrow \Z/4$ be an $i$-invariant quadratic form associated to $e$. The formulas
\[ b(\varepsilon \alpha, \varepsilon \beta )=e(\alpha , \varepsilon \beta )\quad,\quad Q_q(\varepsilon \alpha )=q(\alpha )-Q(\alpha )\qquad \hbox{for }\ \alpha,\beta   \in A_2 \, ,\]
define on $A_i=\varepsilon A_2$ a non-degenerate symmetric form $b$  and a quadratic form $\ Q_q:A_i\rightarrow \Z/4\ $ associated with $ b $. \end{lem}

\pr Since $A_i=\Ker \varepsilon $ is isotropic for $e$, the expression $e(\alpha , \varepsilon \beta )$ is a bilinear function $b$ of $\varepsilon \alpha $ and $\varepsilon \beta $; it is symmetric by $b)$. If $e(\alpha , \varepsilon \beta )=0$ for all $\beta $ in $A_2$ we have $\alpha \in A_i$ because $A_i$ is Lagrangian, hence $\varepsilon \alpha =0$, so $b$ is non-degenerate.

\par  Put $\tilde Q_q(\alpha )=q(\alpha )-Q(\alpha )\in \Z/4$ for $\alpha \in A_2$. We have
\[\tilde Q_q(\alpha +\beta )= \tilde Q_q(\alpha ) + \tilde Q_q(\beta ) +2e( \alpha ,\varepsilon \beta )\ .\]
Take $\beta =\varepsilon \gamma $. Since $q$ is $i$-invariant we have $q(\varepsilon \gamma )=2e(\gamma ,i\gamma )=Q(\varepsilon \gamma )$ by $c)$, hence $\tilde Q_q(\varepsilon \gamma )=0$ and $\tilde Q_q(\alpha +\varepsilon \gamma )=\tilde Q_q(\alpha )$. Thus $\tilde Q_q$  defines a quadratic form $Q_q$ on $A_i$ associated to $b$.\qed

\medskip	
\par  Let $\mathcal{Q}_e^{(i)}$ be the set of  $i$-invariants quadratic forms on $A_2$ associated to $e$. If $q\in\mathcal{Q}_e^{(i)}$ and $\alpha \in A_2$, we have $\alpha +q\in \mathcal{Q}_e^{(i)}$ if and only if $\alpha $ belongs to $A_i^\perp=A_i$; in other words, $\mathcal{Q}_e^{(i)}$ is an affine subspace of $\Q_e$, with direction $A_i$. 

\begin{lem}
The map $q\mapsto Q_q$ is an affine isomorphism of $\mathcal{Q}_e^{(i)}$ onto $\mathcal{Q}_b$.
\end{lem}\label{q}
\pr  We just have to prove the equality $Q_{\alpha +q}=\alpha +Q_q$ for $q\in \mathcal{Q}_e^{(i)}$, $\alpha \in A_i$. Let $\beta \in A_i$;
we write  $\beta  =\varepsilon \beta'$ for some $\beta '\in A_2$. Then
\[Q_{\alpha +q}( \beta ) = 2e(\alpha ,\beta ') + q(\beta ')-Q(\beta ')=2b(\alpha ,\beta )+Q_q(\beta )\ .\qed\]

\begin{rem}\label{syz}
Let $\alpha \in A_2$; we have $b(\varepsilon \alpha ,\varepsilon \alpha )= e(\alpha ,\varepsilon \alpha )=e(\alpha ,i\alpha )\equiv Q(\alpha )\ (\mathrm{mod.}\ 2)$, hence \emph{the form $b$ is symplectic if and only if $\Gamma $ is even}. In this case we have $e(\alpha ,i\alpha )=0$ for all $\alpha \in A_2$; it follows that $\mathcal{Q}_e^{(i)}$ is the set of forms vanishing on $A_i$. Since $A_i$ is  Lagrangian for $e$, this implies that these forms, viewed as quadratic forms
 $A_2\rightarrow \Z/2$, are all  even (that is, their Arf invariant is $0$).
\end{rem}

\medskip	

\section{$i$-invariant theta divisors}
\subsection{Reminder on theta characteristics}\label{theta}
We first recall the classical theory of theta characteristics on an arbitrary p.p.a.v.
$A  =V/\Gamma $. Let $A_2\cong \Gamma /2\Gamma $ be the 2-torsion subgroup of $A$,
 $\mathcal{T}$  the set of symmetric theta divisors on $A$, and $\mathcal{Q}_e$ the set of quadratic forms on $A_2$ associated to the Weyl pairing $e$. The $\F_2$-vector space $A_2$ acts on $\mathcal{T}$ by translation, and on $\mathcal{Q}_e$ by the action defined in (\ref{f2}); both sets are affine spaces over $A_2$, and there is a canonical affine isomorphism $q\mapsto \Theta _q$ of $\mathcal{Q}_e$ onto $\mathcal{T}$. It  can be defined as follows (\cite{M}, \S 2). Let $\gamma \in \Gamma $, and let $\bar\gamma $ be its class in $A_2$. For  $z\in V$, we put 
\[e_\gamma (z)=i^{q(\bar\gamma )}e^{ \pi H(\gamma ,z+\frac{\gamma }{2})}\ .\]
We define an action of $\Gamma $ on the trivial bundle $V\times \C$ by $\gamma .(z,t)=(z+\gamma ,e_\gamma (z)t)$; then the quotient of $V\times \C$ by this action is the line bundle $\O_A(\Theta _q)$ on $A$. 

\subsection{The main results}
We go back to the abelian variety $A_\Gamma $ associated to a Gaussian lattice $\Gamma $. We assume that $\Gamma $ is  unimodular.  We use the notation of (\ref{f2i}). The isomorphism $\mathcal{Q}_e\iso\mathcal{T}$ is compatible with the action of $i$, so $i$-invariant theta divisors correspond to forms $q\in\mathcal{Q}_e^{(i)}$.

\par Let $q\in \mathcal{Q}_e^{(i)}$, and let $L$ be the line bundle  $\O_{A_\Gamma }(\Theta _q)$. We have $i^*L\cong L$; we denote by $\iota :i^*L\rightarrow L$ the unique isomorphism  
inducing  the identity of $L_0$. For each $\alpha \in A_i$, $\iota$  induces an isomorphism $\iota(\alpha ):L_\alpha \rightarrow L_\alpha $. 
 \begin{prop}\label{i(a)}
$\iota (\alpha )$ is the homothety of ratio $i^{Q_q(\alpha )}$.
\end{prop}
\pr The isomorphism $\iota^{-1}:L\iso i^*L$ corresponds to a linear automorphism $j$ of $L$ above $i$:
\[\xymatrix{L \ar[r]^{ j }\ar[d] &L \ar[d] \\
A_\Gamma  \ar[r]^{i}  &\ A_\Gamma \ . }\] 
Consider the automorphism $\tilde j:(z,t)\mapsto (iz,t)$ of $\Gamma _\R\times \C$. Since $ e_{i\gamma}(iz) =e_{\gamma}(z)  $, we have $\tilde j (\gamma .(z,t))=(i\gamma).\tilde j(z,t)  $. Thus $ \tilde j$ factors through an isomorphism 
$ L\rightarrow L $ above $i$ which is the identity on $ L_{0} $, hence equal to $ j $; that is, we have a commutative diagram:
\[\xymatrix{\Gamma _\R\times \C \ar[r]^{\tilde j}\ar[d]_\pi &\Gamma _\R\times \C \ar[d]^\pi \\
L \ar[r]^j  & L
}\]
where $ \pi $ is the quotient map.
\par Let $ \alpha\in A_i $, and let $\gamma $ be an element of $ \Gamma $ whose class $(\mathrm{mod.}\  2\Gamma) $ is $ \alpha $. Then $ \delta:=\frac{i\gamma}{2} -\frac{\gamma}{2} $ belongs to $ \Gamma $. We have
\[j(\pi (\frac{\gamma}{2},t))= \pi(\frac{i\gamma}{2},t) =\pi(\frac{\gamma}{2},e_\delta(\frac{\gamma}{2})_{}^{-1} t) \ ,  \]
hence $ \iota(\alpha) =j(\alpha )^{-1} $ is the homothety of ratio $e_\delta(\frac{\gamma}{2}) $. Let $ \beta $ be the class of $ \delta $ in $ A_2 $. Since $ \gamma=-(1+i)\delta $, we have $ \alpha=\varepsilon\beta $, hence
\[e_\delta(\frac{\gamma}{2})=  i^{q(\beta )}e^{  \frac{\pi }{2}  H(\delta ,\gamma +\delta )}
= i^{q(\beta )-H(\delta,\delta)} =  i^{Q_q(\alpha )}\ .\qed \]

\bigskip	
From $\iota :i^*L\rightarrow L$ we deduce an isomorphism $\iota^\flat :L\iso i_*L$, inducing on global sections
 an automorphism  of $H^0(A_\Gamma ,L)$. 
\begin{prop}\label{arf}
$\iota ^\flat$ acts on $H^0(A_\Gamma ,L)$ by multiplication by $e^{\frac{i\pi }{4}(\sigma (Q_q)+g)}$. 
  
\end{prop}

Note that $\sigma (Q_q)\equiv g\ (\mathrm{mod.}\ 2)$ (\cite{B}, Thm. 1.20, (vi)), so this number is a power of $i$. 

\medskip	
\pr Since $\dim H^0(A_\Gamma ,L)=1$ it suffices to compute $\Tr \iota ^\flat$. This is given by  the holomorphic Lefschetz formula \cite{A-B} applied to $(i,\iota )$. Since $H^i(A_\Gamma ,L)=0$ for $i>0$, we find
\[\Tr \iota^\flat =\sum_{\alpha \in A_i}\frac{\Tr \iota (\alpha )}{ (1-i)^g} = (1-i)^{-g} \sum_{\alpha \in A_i} i^{Q_q(\alpha )}\, .\]

We have $ (1-i)^{-g}=2^{-\frac{g}{2}}e^{\frac{i\pi g}{4}}  $ and $ \sum_{\alpha \in A_i} i^{Q_q(\alpha )}=2^{\frac{g}{2}}e^{\frac{i\pi }{4}\sigma(Q_q) } $ (\cite{B}, Thm. 1.20, (xi)), hence the result.\qed

\begin{prop}
Let $\alpha \in A_i$, and let $m_\alpha (\Theta_q)$ be the multiplicity  of $\Theta_q $ at $\alpha $. We have 
\[2m_\alpha (\Theta _q) \equiv \sigma (Q_q)+g- 2Q_q(\alpha )\quad \mathrm{(mod.\ 8)}\, .\]
\end{prop}
\pr Let $\theta $  be a nonzero section of $H^0(A_\Gamma ,L)$.  Choose a local non-vanishing section $s$ of $L$ around $\alpha $. We can write $\theta =fs$ in a neighborhood of $\alpha $, with $f\in \O_{A_\Gamma ,\alpha }$. We have $\iota ^\flat(\theta )=i^k \theta $ with
 $2k  \equiv \sigma (Q_q)+g\ \mathrm{(mod.}\ 8)$ (Proposition \ref{arf}), hence 
 \[(i^*f )\iota ^\flat(s) =i^k  fs\ .\]
We look at this equality in $\mathfrak{m}_\alpha^m L/\mathfrak{m}_\alpha^{m+1}L$, where 
  $\mathfrak{m}_\alpha $ is the maximal ideal of  $\O_{A_\Gamma ,\alpha }$ and $m:=\allowbreak m_\alpha (\Theta )$.
 We have $i^*f=i^mf\  (\mathrm{mod.}\ \mathfrak{m}_\alpha ^{m+1})$, and $\iota ^\flat(s) =\iota (\alpha )s \  (\mathrm{mod.}\ \mathfrak{m}_\alpha L)$. We obtain $i^m \iota (\alpha )=i^k $, hence the result in view of Proposition \ref{i(a)}.\qed

\begin{cor*}
 The number of $i$-invariant theta divisors $\Theta $ with 
$m_0(\Theta )\equiv 2\ \mathrm{(mod.}\ 4)$ is
\[2^{\frac{g}{2} -1}(2^{\frac{g}{2} }-(-1)^{\frac{g}{4} })\quad \hbox{if $\Gamma $ is even, \ and}\quad 2^{g-2}-2^{\frac{g}{2}-1}\cos\frac{\pi g}{4} \quad  \hbox{if $\Gamma $ is odd;}\]
each of these divisors corresponds to a vanishing thetanull.
\end{cor*}

\pr According to the Proposition, we have $m_0(\Theta_q )\equiv 2$ $\mathrm{(mod.}\ 4)$ if and only if  $\sigma (Q_q)\equiv 4-g$ \allowbreak$\mathrm{(mod.}\ 8)$.
When $q$ runs over $\mathcal{Q}_e^{(i)}$, $Q_q$ runs over $\mathcal{Q}_b$ (Lemma \ref{q}),
so we must find how many elements $Q$ of  $\mathcal{Q}_b$ satisfy $\sigma (Q)\equiv 4-g$ \allowbreak$\mathrm{(mod.}\ 8)$.

\par If $\Gamma $ is even (so that $g$ is divisible by $4$), we identify $\mathcal{Q}_b$ with the set of
 quadratic forms $Q: A_2\rightarrow \Z/2$ associated with the symplectic form $b$; the previous congruence 
becomes 
$\mathrm{Arf}(Q)\equiv 1+\frac{g}{4}\ (\mathrm{mod.}\ 2) $. 
There are 
$2^{\frac{g}{2} -1}(2^{\frac{g}{2} }+1)$ such forms with Arf invariant $0$ and $2^{\frac{g}{2} -1}(2^{\frac{g}{2} }-1)$ with Arf invariant 1, hence the result.

Assume that $\Gamma $ is odd; we choose an orthonormal basis $(e_1,\ldots ,e_g)$ for $b$. The forms $Q\in \mathcal{Q}_b$ are determined by their values $Q(e_i)=\pm 1$; the condition is that the number $g^+$ of $+1$ values satisfies 
\[2g^+ -g\equiv 4-g  \ (\mathrm{mod.}\ 8)\ ,\ \hbox{hence }\ g^+\equiv 2\ (\mathrm{mod.}\ 4)\ .\]
The number of forms with the required property is thus the number of subsets  $E\subset \{1,\ldots ,g\}$ with $\mathrm{Card}(E)\equiv 2\ (\mathrm{mod.}\ 4)$, that is
 \[{g\choose 2}+{g\choose 6}+\ldots  = \frac{1}{4}\bigl[(1+1)^g+(1-1)^g-(1+i)^g-(1-i)^g\bigr] =2^{g-2}-2^{\frac{g}{2}-1}\cos\frac{\pi g}{4}\ .  \qed \]

\bigskip	
Thus we find a number of vanishing thetanulls asymptotically equivalent to $2^{g-1}$ when $\Gamma $ is even, and $2^{g-2}$ when $\Gamma $ is odd. These numbers are  rather modest, at least by comparison with the number of vanishing thetanulls of a hyperelliptic Jacobian, which is asymptotically equivalent to $2^{2g-1}$.
However, when $\Gamma $ is even, the vanishing thetanulls of $A_\Gamma $ have the particular property of being ``syzygetic'' in the classical terminology, which just means that the corresponding quadratic forms (\ref{theta}) lie in an affine subspace of $\Q_e$ which consists of even forms (Remark \ref{syz}). Such a subspace has dimension $\leq g$, and it might be that 
the number given by the Corollary in the even case is the maximum possible for a syzygetic subset of vanishing thetanulls. 

\medskip	

\section{Complements}
\subsection{Automorphisms}
The automorphism group of $A_\Gamma $ is the centralizer of $i$ in $\Aut(\Gamma )$. This group can be rather large: it has order 46080 for $\Gamma =E_8$ and 2012774400 for $\Gamma =\Lambda _{24}$ \cite{C-S}. For the lattice $\Gamma _{2g}$ (Example \ref{examples}.1) with $g>4$, it has order $2^{2g-1}g!$.
\par For the lattice $\Gamma =\Gamma _0\otimes_{\Z}\Z[i]$ of Example \ref{examples}.3, $\Aut(A_\Gamma )$ is generated by $i$ and the group $\Aut(\Gamma _0)$. 
 Note that there are examples of  unimodular lattices (even or odd) $\Gamma_0 $ with  $\Aut(\Gamma_0 ) =  \allowbreak \{\pm 1\}$ \cite{Ba}, so that $\Aut(A_\Gamma )$ is reduced to $\{\pm 1,\pm i\}$.

\subsection{Jacobians}
We observe that for $g>1$ the p.p.a.v. $A_\Gamma $ can \emph{not} be a Jacobian. Indeed, let $C$ be a curve of genus $g$; if $JC\cong A_\Gamma $, Torelli theorem provides an automorphism $u$ of $C$ inducing either $ i$ or $-i$ on $JC$, hence also on $T_0(JC)=H^0(C,K_C)^*$. Then $u$ acts trivially on the image of the canonical map $C\rightarrow \P(H^0(C,K_C)^*)$; this implies that $u$ is the identity or that $C$ is  hyperelliptic and  $u$ is the hyperelliptic involution. But in these cases $u$ acts on $H^0(C,K_C)$ by multiplication by $\pm 1$, a  contradiction.\qed

\end{document}